\documentclass[11pt]{article}
\usepackage{amsmath}
\usepackage{amsthm}
\usepackage{amssymb}
\usepackage{mathrsfs}
\usepackage{amssymb,amsmath}
\usepackage{leftidx}
\usepackage{graphicx}
\usepackage{amsfonts}
\usepackage{graphicx, color}
\usepackage{cite}
\usepackage{enumerate}%ÁÐ±íºê°ü
\usepackage{enumitem}
\usepackage{subfigure}
\usepackage{multicol}
\usepackage{multirow}%ºÏ²¢¶àÐÐµ¥Ôª¸ñ
\usepackage{booktabs}%´´½¨Ã»ÓÐÊúÏß·Ö¸îµÄ±í¸ñ
\usepackage{algpseudocode,algorithm,algorithmicx}%Éú³ÉÎ±´úÂëËã·¨Í¼
\newtheorem{theorem}{Theorem}[section]
\newtheorem{lemma}{Lemma}[section]

\newtheorem{example}{Example}

\newcommand{\zd}{\,\mathrm{d}}

\newcommand{\diff}{\triangledown_{\tau}}
\newcommand{\abs}[1]{\left|#1\right|}
\newcommand{\absb}[1]{\big|#1\big|}

\newcommand{\bra}[1]{\left(#1\right)}
\newcommand{\brab}[1]{\big(#1\big)}
\newcommand{\braB}[1]{\Big(#1\Big)}
\newcommand{\brat}[1]{(#1)}
\newcommand{\kbra}[1]{\left[#1\right]}
\newcommand{\kbrab}[1]{\big[#1\big]}

\newcommand{\myinner}[1]{\left\langle#1\right\rangle}
\newcommand{\myinnerb}[1]{\big\langle#1\big\rangle}

\newcommand{\mynorm}[1]{\left\|#1\right\|}
\newcommand{\mynormb}[1]{\big\|#1\big\|}

\begin{document}
\title{The variable-step L1 scheme preserving a compatible energy law
for time-fractional Allen-Cahn equation}
\author{
Hong-lin Liao\thanks{ORCID 0000-0003-0777-6832; Department of Mathematics,
Nanjing University of Aeronautics and Astronautics,
Nanjing 211106, P.R. China. Hong-lin Liao (liaohl@csrc.ac.cn and liaohl@nuaa.edu.cn)
is supported by a grant 12071216 from National Natural Science Foundation of China.}
\quad Xiaohan Zhu\thanks{Department of Mathematics, Nanjing University of Aeronautics and Astronautics,
211106, P.R. China. Email: cyzhuxiaohan@163.com.}
\quad Jindi Wang\thanks{School of Mathematics and Computational Science,
Xiangtan University, Xiangtan 411105, P.R. China.
Jindi Wang (wangjindixy@163.com) is supported by a grant
XDCX2020B078 from Hunan Provincial Innovation Foundation for Postgraduate.}
}
%%%%%%%%%%%%%%%%%%%%%%%%%%%%%%%%%%%%%%%%%%%%%%%%%%%%%%%%%%%%%%%%%%%%%%%%%%%%%%%%%%%%%%%%%%
\date{December 24, 2020}
\maketitle
\normalsize

\begin{abstract}
  In this work, we revisit the adaptive L1 time-stepping scheme for solving
  the time-fractional Allen-Cahn equation in the Caputo's form.
  The L1 implicit scheme is shown to preserve a variational energy dissipation law
  on arbitrary nonuniform time meshes by using the recent discrete analysis tools, i.e.,
  the discrete orthogonal convolution kernels and discrete complementary convolution kernels.
  Then the discrete embedding techniques and the fractional Gr\"onwall
  inequality were applied to establish an $L^2$ norm error estimate on nonuniform time meshes.
  An adaptive time-stepping strategy according to the dynamical feature of the system
  is presented to capture the multi-scale behaviors
  and to improve the computational performance.
  \noindent{\bf AMS subject classiffications.}\;\; 35Q99, 65M06, 65M12, 74A50
\end{abstract}

%%%%%%%%%%%%%%%%%%%%%%%%%%%%%%%%%%%%%%%%%%%%%%%%%%%%%%%%%%%%%%%%%%%%%%%%%%%%%%%%%%%%%%%
\section{Introduction}

We consider the numerical approximations for
time-fractional Allen-Cahn (TFAC) equation
\begin{align}\label{cont:FAC equ}
\partial_{t}^{\alpha}\Phi=-\kappa\mu
\quad\text{where the potential}\;\;\mu:=f(\Phi)-\epsilon^2\Delta\Phi,
\end{align}
on a bounded regular domain
$\mathbf{x}\in\Omega\subseteq\mathbb{R}^2$
subject to periodic boundary conditions.
Here, $\epsilon>0$ is an interface width parameter,
$\kappa>0$ is the mobility coefficient,
and the nonlinear bulk force $ f(\Phi) $ is taken as the polynomial
double-well potential $f(\Phi)=\Phi^3-\Phi$.
The notation $\partial_{t}^{\alpha}:={}_{0}^{C}\!D_{t}^{\alpha}$ in {\eqref{cont:FAC equ}}
represents the fractional Caputo derivative of order $\alpha$
with respect to $t$, that is,
\begin{align}\label{def:Caputo deriva}
(\partial_{t}^{\alpha}v)(t)
:=(\mathcal{I}_{t}^{1-\alpha}v')(t)
\quad \text{for $0<\alpha<1$},
\end{align}
in which the fractional  Riemann-Liouville integral $\mathcal{I}_{t}^{\beta}$
of order $\beta>0$ is given by
\begin{align} \label{RLDef}
(\mathcal{I}_{t}^{\beta}v)(t)
:=\int_{0}^{t}\omega_{\beta}(t-s)v(s)\zd{s}\quad\text{where\; $\omega_{\beta}(t):=t^{\beta-1}/\Gamma(\beta)$.}
\end{align}

As well known, the energy dissipation law is an important
and essential property of the  classical phase field models.
Recall the following Ginzburg-Landau energy functional \cite{Allen1979A},
\begin{align}\label{cont:free energy}
E[\Phi]:=\int_{\Omega}\braB{\frac{\epsilon^2}{2}\abs{\nabla\Phi}^2
+F(\Phi)}\zd\mathbf{x}\quad \text{where\; $F(\Phi)=\frac{1}{4}\bra{\Phi^2-1}^2$.}
\end{align}
The classical Allen-Cahn (AC) model preserves
the following energy dissipation law,
\begin{align}\label{cont:classical energy law}
\frac{\zd E}{\zd t}
+\kappa\mynormb{\tfrac{\delta E}{\delta \Phi}}^2=0\quad \text{for $t>0$},
\end{align}
%and the maximum bound principle,
%\begin{align*}%\label{FracMaximumPrinciple}
%|\Phi(\mathbf{x},t)|\le{1}\;\text{for $t>0$}\quad
%\text{if}\quad|\Phi(\mathbf{x},0)|\le{1}.
%\end{align*}
where the inner product $\bra{u,v}:=\int_{\Omega}uv\zd\mathbf{x}$,
and the associated $L^2$ norm $\mynormb{u}:=\sqrt{(u,u)}$
for all $u,v\in L^2(\Omega).$
It is of great interest to design some numerical algorithms
that preserve the energy dissipation law at each time level
because non-energy-stable numerical schemes would not accurately capture the coarsening dynamics
or lead to numerical instability.

For the classical gradient flows, there are several effective
strategies to develop energy stable numerical algorithms,
 such as the convex splitting method \cite{wise2009an,Cheng2019An},
stabilization technique \cite{Xu2006Stability,Shen2010Numerical},
invariant energy quadratization approach
\cite{Gong2018lLnear,Gong2020Arbitrarily}
and scalar auxiliary variable formulation
\cite{Jie2018The}.
Compared  with the classical phase field models, however, the theoretical works regarding
the energy stable property of the time-fractional phase field models are limited.
%and \cite[Theorem 2.3]{Du2019Time},
It was shown \cite[Theorem 4.2]{Tang2018On} that the TFAC model \eqref{cont:FAC equ}
admits the maximum bound principle
\begin{align}\label{cont: TFAC maximum bound principle}
\absb{\Phi(x,t)}\le 1\quad
\text{if}\quad \absb{\Phi(x,0)}\le 1\quad\text{for $t>0$,}
\end{align}
and preserves the following global  energy stable property
\begin{align}\label{cont:global energy law}
E\kbra{\Phi(t)}\le E\kbra{\Phi(0)}\quad \text{for $t>0$.}
\end{align}
It implies that the energy is bounded by the initial one.
At the discrete levels, the authors in \cite{Tang2018On}
combined the uniform L1 formula with stabilization technique
to develop a numerical scheme preserving this global energy stability.

Very recently, two nonlocal energy decaying laws for the time-fractional
phase field models were developed in \cite{Quan2020How},
including a time-fractional energy dissipation law
\begin{align}\label{cont:frac energy law}
\bra{\partial_t^\alpha E}(t)\le 0\quad \text{for $t>0$},
\end{align}
and a weighted energy dissipation law,
\begin{align}\label{cont:weight energy law}
\frac{\zd E_{\omega}}{\zd t}\le 0\quad \text{for $t>0$}
\quad \text{where\; $E_{\omega}(t):=\int_0^1\omega(\theta)E(\theta t)\zd \theta$},
\end{align}
where $\omega(\theta)\ge0$ is some weight function satisfying
$\int_0^1\omega(\theta)\zd t=1$ and $E(\theta t)=E\kbra{\Phi\bra{\cdot,\theta t}}$
is the classical energy defined by \eqref{cont:free energy}.
The authors also proposed some numerical approaches,
including the convex-splitting and
scalar auxiliary variable schemes, in \cite{Quan2020Numerical}
to preserve the above two energy dissipation laws.
It seems that the time-fractional energy dissipation law \eqref{cont:frac energy law}
is consistent with an energy decaying law, that is, $\bra{\partial_t^\alpha E}(t) \rightarrow
\bra{\partial_t E}(t)\le 0$ as $\alpha\rightarrow 1$;
%\begin{align*}
%\bra{\partial_t^\alpha E}(t) \rightarrow
%\bra{\partial_t E}(t)\le 0 \quad \text{as $\alpha\rightarrow 1$};
%\end{align*}
while the weighted law \eqref{cont:weight energy law}
may be not compatible with the classical one as
$\alpha\rightarrow 1$.

In this  paper, we shall view the TFAC equation \eqref{cont:FAC equ}
as a time-fractional gradient flow since it can recover the classical
AC model when the fractional order $\alpha\rightarrow 1$.
Recently, Liao, Tang and Zhou \cite{Liao2020An} explored a variational energy functional
\begin{align}\label{cont:TFAC variational energy}
\mathcal{E}_\alpha[\Phi]:=E[\Phi]
+\frac{\kappa}{2}\mathcal{I}_{t}^{\alpha}\mynormb{\tfrac{\delta E}{\delta \Phi}}^2,
\end{align}
which was proven in \cite[Section 1.1]{Liao2020An} to satisfy a variational energy dissipation law,
\begin{align}\label{cont:TFAC energy law}
\frac{\zd \mathcal{E}_\alpha}{\zd t}
+\frac{\kappa}{2}\omega_{\alpha}(t)\mynormb{\tfrac{\delta E}{\delta \Phi}}^2\le0
\quad \text{for $t>0$.}
\end{align}
Obviously, it leads to the global energy law \eqref{cont:global energy law} directly.
As remarked in \cite{Liao2020An}, this type of energy law seems naturally because it
asymptotically compatible with the classical energy dissipation law.
Actually,  as the fractional order $\alpha\rightarrow 1$,
the variational energy dissipation law \eqref{cont:TFAC energy law} naturally approaches
the energy dissipation law \eqref{cont:classical energy law},
\begin{align*}
\frac{\zd E}{\zd t}+\kappa\mynormb{\tfrac{\delta E}{\delta \Phi}}^2\le0
\quad \text{for $t>0$}.
\end{align*}
It is worthy mentioning that the new energy law \eqref{cont:TFAC energy law} was derived via the  Riemann-Liouville version
of the TFAC model \eqref{cont:FAC equ},
\begin{align}\label{cont:TFAC equi RL}
\partial_t\Phi=-\kappa{}^{R}\!\partial_t^{1-\alpha}\mu
\quad\text{with}\quad\mu:=\tfrac{\delta E}{\delta \Phi}=f(\Phi)-\epsilon^2\Delta\Phi,
\end{align}
which can be obtained by acting the Riemann-Liouville fractional derivative
${}^{R}\!\partial_t^\beta:=\partial_t\mathcal{I}_{t}^{1-\beta}$
on both sides of the equation \eqref{cont:FAC equ}
and using the semigroup property
$\mathcal{I}_{t}^{\alpha}\mathcal{I}_{t}^{\beta}=\mathcal{I}_{t}^{\alpha+\beta}$.
Also, a nonuniform L1-type (called L1$_R$) time-stepping scheme with the approximation order $1+\alpha$
preserving the maximum bound principle \eqref{cont: TFAC maximum bound principle}
and the new variational energy dissipation law \eqref{cont:TFAC energy law} was investigated in \cite{Liao2020An}
for the Riemann-Liouville version \eqref{cont:TFAC equi RL}.
%However, no convergence results were established.

%For the gradient flow systems involving the multiple time scale
%behaviors such as the TFAC model \eqref{cont:FAC equ},
%adaptive time stepping would be of particular interesting to capture delicate
%small scale structures with correct physical information
%and to accelerate the long-time simulations.
We shall consider a direct approximation of the TFAC model \eqref{cont:FAC equ}
with the  L1 formula of Caputo derivative \eqref{def:Caputo deriva}.
For a finite $T>0$,
consider  $0=t_{0}<t_{1}<\cdots<t_{k}<\cdots<t_{N}=T$.
Let the variable time-steps $\tau_{k}:=t_{k}-t_{k-1}$ for $1\le{k}\le{N}$,
the maximum step size $\tau:=\max_{1\le{k}\le{N}}\tau_{k}$,
and the adjoint time-step ratios $r_k:=\tau_k/\tau_{k-1}$ for $2\le k\le N$.
%$$r_k:=\tau_k/\tau_{k-1}\quad \text{for $2\le k\le N$}.$$
Given a grid function $\{v^{k}\}_{k=0}^N$, let $\diff v^k:=v^k-v^{k-1}$ and
$\partial_\tau v^k:=\diff v^k/\tau_k$ for $k\geq{1}$.
The nonuniform L1 formula of Caputo derivative \eqref{def:Caputo deriva} reads \cite{Liao2018Sharp,Liao2018Unconditional},
\begin{align}\label{sche:L1 formula}
(\partial_{\tau}^{\alpha}v)^n
:=\sum_{k=1}^{n}a_{n-k}^{(n)}\diff v^{k}
\quad \text{with}\quad
a_{n-k}^{(n)}
:=\frac{1}{\tau_{k}}\int_{t_{k-1}}^{t_{k}}\omega_{1-\alpha}(t_n-s)\zd{s},\;\;1\le k\le n.
\end{align}
%for $1\le{k}\le{n}$.
By using finite difference approximation in space (see section 2),
we have the following L1 implicit scheme subject to the initial data
$\phi_h^0=\Phi_0(\mathbf{x})$ and periodic boundary conditions,
\begin{align}\label{sche:adap AC L1}
(\partial_{\tau}^{\alpha}\phi_h)^{n}&=-\kappa\mu_h^n \quad\text{with}\quad
\mu_h^n:=f(\phi_h^n)-\epsilon^2\Delta_h \phi_h^n\quad \text{for $n\ge1$.}
\end{align}

In fact, this backward Euler-type scheme \eqref{sche:adap AC L1} has been investigated in our previous work \cite{Ji2020Simple}.
As noticed, two types of nonuniform L1 schemes,
including the first-order stabilized semi-implicit method and the $(2-\alpha)$-order implicit scheme \eqref{sche:adap AC L1},
for the TFAC model \eqref{cont:FAC equ} have been proven to
preserve the maximum bound principle \eqref{cont: TFAC maximum bound principle},
see \cite[Theorem 2.1 and Theorem 2.2]{Ji2020Simple}.
Nonetheless, no any discrete energy dissipation laws were established
on nonuniform meshes.

By using the positive definiteness of L1 kernels with the uniform time step,
essentially due to the key result \cite[Proposition 5.2]{Lopez1990A},
the first-order stabilized semi-implicit scheme was proved in \cite[Lemma 2.6]{Ji2020Simple}
to preserve the global energy law \eqref{cont:global energy law}.
It is to mention that,  by using the main theorem \cite[Theorem 1.1]{Liao2020Positive} on
the positive definiteness of real quadratic form with
variable coefficients, the first-order stabilized scheme was shown \cite[Proposition 4.2]{Liao2020Positive} to preserve
the global energy dissipation law \eqref{cont:global energy law} on arbitrary time meshes.
However, no any discrete energy dissipation laws have been established for the implicit scheme
\eqref{sche:adap AC L1}, even on the uniform grid.

This paper aims to fill this gap for the backward Euler-type scheme
\eqref{sche:adap AC L1} with variable time-steps
by establishing a discrete energy dissipation law that is
asymptotically compatible with the discrete energy law
of the backward Euler scheme (cf. \cite[(2.3)]{xu2019on}) for the AC model,
\begin{align}\label{sche: backward Euler AC}
\partial_\tau\phi_h^n=-\kappa\mu_h^n\quad\text{with}\quad
\mu_h^n:=f\brab{\phi_h^n}-\epsilon^2\Delta_h\phi_h^n
\quad \text{for $n\ge1$.}
\end{align}
For this simple scheme, \cite[Theorem 2.1]{xu2019on} stated (by our notations, such as $\kappa=1/\epsilon^2$) that
\begin{align}\label{cond: time-step size}
\text{if the step size $\tau_n\le 1/\kappa$, the backward Euler scheme is convex and uniquely solvable,}
\end{align}
and satisfies the following energy dissipation law
\begin{align}\label{ieq: energy dissipation law}
\partial_\tau E\kbra{\phi^n}
+\frac{\kappa}{2}\mynormb{\mu^n}^2\le 0\quad \text{for $n\ge1$.}
\end{align}
In this sense, for the TFAC model \eqref{cont:FAC equ} with Caputo fractional derivative,
it would be the first work on the direct approximation with variable time-steps that can preserve both
the maximum bound principle and the energy dissipation law at each time level.

Our main tool for constructing the discrete energy dissipation law is the so-called discrete orthogonal convolution (DOC)
kernels $ {\theta_{n-k}^{(n)}} $ defined by the following  recursive procedure
\begin{align}\label{def:DOC kernel}
{\theta_{0}^{(n)}}:=\frac{1}{a_{0}^{(n)}}
\quad \mathrm{and} \quad
{\theta_{n-k}^{(n)}}:=-\frac{1}{a_{0}^{(k)}}
\sum_{j=k+1}^n {\theta_{n-j}^{(n)}}a_{j-k}^{(j)}
\quad \text{for $1\le k\le n-1$}.
\end{align}
Obviously, they satisfy the following discrete orthogonal identity
\begin{align}\label{DOC:orthogonal identity}
\sum_{j=k}^n {\theta_{n-j}^{(n)}}a_{j-k}^{(j)}
\equiv \delta_{nk}\quad\text{for $1\le k\le n$,}
\end{align}
where $\delta_{nk}$ is the Kronecker delta symbol. By acting the DOC kernels
on the L1 formula \eqref{sche:L1 formula} and applying the discrete orthogonal identity
\eqref{DOC:orthogonal identity}, one gets
\begin{align}\label{sche: DOC action L1 formula}
\sum_{j=1}^{n}\theta_{n-j}^{(n)}(\partial_{\tau}^{\alpha}v)^j
=&\,\sum_{k=1}^{n}\diff v^k\sum_{j=k}^n {\theta_{n-j}^{(n)}}a_{j-k}^{(j)}
\equiv\diff v^n\quad\text{for $1\le n\le N$.}
\end{align}
Then, by acting the DOC kernels
$\theta_{n-j}^{(n)}$ on the L1 scheme \eqref{sche:adap AC L1}, this identity introduces
the following equivalent scheme with respect to the DOC kernels,
\begin{align}\label{sche:adap AC L1 equi form}
\diff\phi_h^n=-\kappa\sum_{j=1}^n\theta_{n-j}^{(n)}\mu_h^j
\quad\text{with}\quad
\mu_h^n:=\bra{\phi_h^n}^3-\phi_h^n-\epsilon^2\Delta_h\phi_h^n
\quad\text{for $1\le n\le N$.}
\end{align}
Actually, the original L1 scheme \eqref{sche:adap AC L1} can be recovered
from the equivalent formulation \eqref{sche:adap AC L1 equi form} by acting $a_{m-n}^{(m)}$
on both sides of \eqref{sche:adap AC L1 equi form}
and using the following mutually orthogonal identity (by the proof of \cite[Lemma 2.1]{Liao2020Positive})
 \begin{align}\label{eq: mutual orthogonal identity}
 \sum_{j=k}^na_{n-j}^{(n)}\theta_{j-k}^{(j)}\equiv \delta_{nk}\quad\text{for $1\le k\le n$.}
   \end{align}

We will use the formulation \eqref{sche:adap AC L1 equi form} to derive
the desired discrete energy law as
the derivation of its continuous counterpart \eqref{cont:TFAC energy law} in \cite[section 1.1]{Liao2020An}.
Obviously, this formulation can also be viewed
as a direct numerical approximation for the
Riemann-Liouville version \eqref{cont:TFAC equi RL} of the TFAC model \eqref{cont:FAC equ}.
In this sense, the DOC kernels $\theta_{n-k}^{(n)}$ ``define"
an indirect formula of the Riemann-Liouville fractional derivative ${}^{R}\!\partial_t^{1-\alpha}v$ by
\begin{align}\label{sche: indirect L1 formula RL}
\bra{{}^{R}\!\partial_\tau^{1-\alpha}v}^n
:=\frac{1}{\tau_n}\sum_{j=1}^n\theta_{n-j}^{(n)}v^j
\quad \text{for $n\ge1.$}
\end{align}
With the help of the mutually orthogonal identity \eqref{eq: mutual orthogonal identity},
one can recover the L1 formula \eqref{sche:L1 formula} by acting the L1 kernels $a_{m-n}^{(m)}$
on both sides of \eqref{sche: DOC action L1 formula}. That is to say, the DOC kernels $\theta_{n-k}^{(n)}$ also define a
\emph{reversible discrete transformation} between
the nonuniform L1 formula \eqref{sche:L1 formula} of
Caputo derivative $\partial_{t}^{\alpha}$ and
the indirect formula \eqref{sche: indirect L1 formula RL}
of  Riemann-Liouville  derivative ${}^{R}\!\partial_t^{1-\alpha}$.

In the next section,
the unique solvability of the suggested nonuniform L1 method is proved in Theorem \ref{thm:uni solva}.
Theorem \ref{thm:energy decay law} establishes the discrete variational energy dissipation law
for the L1 scheme by using the DOC kernels \eqref{def:DOC kernel} and the
so-call discrete complementary convolution (DCC) kernels.
It is to emphasize that the unique solvability is proved
by using the original numerical scheme \eqref{sche:adap AC L1},
while, as mentioned, the energy stability is established by using
the equivalent convolution form \eqref{sche:adap AC L1 equi form}.

By making use of the discrete $H^1$ norm solution bound obtained from
the discrete energy stability,
an $L^2$ norm error estimate is then achieved in section 3
with the help of the discrete fractional Gr\"onwall
inequality. Numerical examples including the accuracy verification
and simulations of coarsening dynamics are carried out in section 4
to illustrate the effectiveness of the L1 scheme,
especially when it is coupled with an adaptive time-stepping strategy.

%In summary, our contributions in this paper are three  folds:
%\begin{itemize}
%  \item By virtue of  the recent DOC and DCC techniques,
%   the fully implicit adaptive time-stepping L1 scheme is prove to preserve
%   the variational energy dissipation property \eqref{cont:TFAC energy law}
%   in the discrete levels.
% \item  The sharp $L^2$ error estimate is established
%  on a general class of nonuniform time meshes so that the adaptive time steps
%  can be used to capture the multi-scale dynamic behaviors.
%  \item Several numerical examples are presented to confirmed our theoretical findings.
%\end{itemize}

Throughout this paper, any subscripted $C$,
such as $C_v$, $C_\phi$ and $C_{\gamma}$,  denotes a generic positive constant,
not necessarily the same at different occurrences;
while, any subscripted $c$, such as $c_\Omega,c_0,c_1,c_2$ and $c_3$,
denotes a fixed positive constant. Always, the appeared constants are  dependent on the given data
and the solution, but always independent of the spatial lengths, the time $t_n$,
the time-step sizes $\tau_n$ and time-step ratios $r_n$.

\section{Solvability and energy dissipation law}
For simplicity,
we cover $\Omega=(0,L)^2$ by
the discrete grid
$\bar{\Omega}_h:=\big\{\mathbf{x}_{h}=(ih,jh)\,|\,0\le i,j\le M\big\}$
with the uniform length $h:=L/M$ for some integer $M$.
Let $\Omega_h:=\bar{\Omega}_h\cap\Omega$ and
denote the space of $L$-periodic grid functions
$\mathbb{V}_{h}:=\{v\,|\,v=\bra{v_h}\; \text{is $L$-periodic for}\; \mathbf{x}_h\in\bar{\Omega}_h\}.$

For any grid functions $v,w\in\mathbb{V}_h$,
define the discrete inner product
$\myinner{v,w}:=h^2\sum_{\mathbf{x}_h\in\Omega_{h}}v_h w_h$,
the associated $L^{2}$ norm $\mynorm{v}:=\sqrt{\myinner{v,v}}$
and $L^p$ norm $\mynorm{v}_{\ell^p}=\sqrt[p]{h^2\sum_{\mathbf{x}_h\in\Omega_{h}}|v_h|^p}$.
The standard second-order central finite difference approximations
are used in the space discretization.
Let $\nabla_h$ and $\Delta_h$ be the discrete gradient and Laplace
operations  in the point-wise sense such that the discrete Green's formula
 $\myinner{-\Delta_hv,w}=\myinner{\nabla_hv,\nabla_hw}$ holds.

\subsection{Unique solvability}

%%%%%%%%%%%%%%%%%%%%%%%%%%%%%%%%%%%%%%%%%%%%%%%%%%%%%%%%%%%%%%
\begin{theorem}\label{thm:uni solva}
Under the time-step restriction
\begin{align}\label{restr:solvab time step}
\tau_n\le\frac{1}{\sqrt[\alpha ]{\kappa\Gamma(2-\alpha)}}	
\end{align}
the fully implicit L1 scheme \eqref{sche:adap AC L1} is uniquely solvable.
\end{theorem}
%%%%%%%%%%%%%%%%%%%%%%%%%%%%%%%%%%%%%%%%%%%%%%%%%%%%%%%%%%%%%%

\begin{proof}
Consider the following  energy functional $G[z]$,
\begin{align*}
G[z]:=\frac{a_0^{(n)}}{2}\mynorm{z-\phi^{n-1}}^2
+\myinnerb{\mathcal{L}^{n-1},z}
+\frac{\epsilon^2}{2}\kappa \mynormb{\nabla_h z}^2+\frac{\kappa}{4} \mynorm{z}_{\ell^4}^4-\frac{\kappa}{2}\mynormb{z}^2\quad\text{for $n\ge1$,}
\end{align*}
where we denote $\mathcal{L}^{n-1}:=\sum_{k=1}^{n-1}a_{n-k}^{(n)}\diff\phi^k$.
The solution of nonlinear equation \eqref{sche:adap AC L1}
is equivalent to the minimum of $G[z]$ if and only if it is strictly convex
and coercive on $\mathbb{V}_h$, see \cite{Feng2012Analysis}.

In details, the time-step restriction condition \eqref{restr:solvab time step}
shows that $a_0^{(n)} \geq \kappa$.
So it can be easily verified that the functional $G[z]$
is convex with respect to $z$ on $\mathbb{V}_h$,
\begin{align*}
\frac{\zd^2G}{\zd s^2}[z+s\psi]\Big|_{s=0}
&=a_0^{(n)}\mynormb{\psi}^2+\kappa\epsilon^2\mynormb{\nabla_h\psi}^2
+3\kappa\mynormb{z\psi}^2-\kappa\mynormb{ \psi}^2 \\
&=(a_0^{(n)}-\kappa)\mynormb{ \psi }^2
+\kappa\epsilon^2\mynormb{\nabla_h\psi}^2+3\kappa\mynormb{z\psi}^2>0.
\end{align*}
Moreover, it is straightforward to show
that the functional $G[z]$ is coercive on $\mathbb{V}_h$, that is,
\begin{align*}
G[z]&\ge \myinnerb{\mathcal{L}^{n-1},z}
+\frac{\kappa}{4}\mynorm{z}_{\ell^4}^4-\frac{\kappa}{2}\mynormb{z}^2
\ge \frac{\kappa}{4}\mynorm{z}_{\ell^4}^4-\kappa\mynormb{z}^2
-\frac{1}{2\kappa}\mynormb{\mathcal{L}^{n-1}}^2\\
&\ge\frac{\kappa}{2}\mynormb{z}^2-\frac{9\kappa}{4}\abs{\Omega}
-\frac{1}{2\kappa}\mynormb{\mathcal{L}^{n-1}}^2,
\end{align*}
where the inequality
$\mynorm{v}_{\ell^4}^4\ge 6\mynormb{v}^2-9\abs{\Omega}$
has been used in the last step.
Hence, the functional $G[z]$  has a unique minimizer
and then the scheme \eqref{sche:adap AC L1} is uniquely solvable.
%It completes the proof.
\end{proof}

Let the fractional order $\alpha\rightarrow 1$, the L1 kernels
$a_0^{(n)}\rightarrow 1/\tau_n$ and $a_{n-k}^{(n)}\rightarrow 0$ for $1\le k\le n-1$.
The L1 scheme \eqref{sche:adap AC L1} degrades into the backward Euler scheme \eqref{sche: backward Euler AC}.
Note that, the time-step restriction \eqref{restr:solvab time step}
for the unique solvability approaches $\tau_n\le 1/\kappa$,
just the time-step condition \eqref{cond: time-step size}
for the unique solvability of \eqref{sche: backward Euler AC}.

\subsection{Discrete energy dissipation law}

We have the following result on the L1 kernels $a_{n-k}^{(n)}$ defined in \eqref{sche:L1 formula}.

%%%%%%%%%%%%%%%%%%%%%%%%%%%%%%%%%%%%%%%%%%%%%%%%%%%%%%%%%%%%%%
\begin{lemma}\cite[Proposition 4.1]{Liao2020Positive}\label{lem:L1 kernel property}
For $n\ge 2$, the L1 kernels $a_{j}^{(n)}$ in \eqref{sche:L1 formula} satisfy
\begin{enumerate}%[(i)]
\item[(i)]    $a_{j-1}^{(n)}>a_{j}^{(n)}>0$ for $1\le j\le n-1$;
  \item[(ii)]    $a_{j-1}^{(n-1)}>a_{j}^{(n)}$ for $1\le j\le n-1$;
  \item[(iii)]   $a_{j-1}^{(n-1)}a_{j+1}^{(n)}>a_{j}^{(n-1)}a_{j}^{(n)}$ for $1\le j\le n-2$.
\end{enumerate}
\end{lemma}
%%%%%%%%%%%%%%%%%%%%%%%%%%%%%%%%%%%%%%%%%%%%%%%%%%%%%%%%%%%%%%

By using Lemma \ref{lem:L1 kernel property}, we can follow the proof of
\cite[Lemma 2.3]{Liao2020Positive} to give the following result on the associated DOC kernels $ {\theta_{n-k}^{(n)}}$.

%%%%%%%%%%%%%%%%%%%%%%%%%%%%%%%%%%%%%%%%%%%%%%%%%%%%%%%%%%%%%%%%%%%%%%%%%%%%
\begin{lemma}\label{lem:DOC kernel property}
For any $n\geq 2$, the DOC kernels ${\theta_{n-k}^{(n)}}$ defined in
\eqref{def:DOC kernel} satisfy
\begin{align*}
\theta_{0}^{(n)}>0\quad\text{and}\quad \theta_{n-k}^{(n)}<0\;\;\text{ for $1\le k\le n-1$};\quad
\text{but}\quad \sum_{k=1}^n {\theta_{n-k}^{(n)}}>0.
\end{align*}
%\begin{itemize}
%  \item [(i)] $\displaystyle \theta_{0}^{(n)}>0$ and $ {\theta_{n-k}^{(n)}}<0$ for $1\le k\le n-1$;
%\item [(ii)] $ \displaystyle \sum_{k=1}^n {\theta_{n-k}^{(n)}}>0.$
%\end{itemize}
\end{lemma}
%%%%%%%%%%%%%%%%%%%%%%%%%%%%%%%%%%%%%%%%%%%%%%%%%%%%%%%%%%%%%%%%%%%%%%%%%%%%

To derive the discrete energy law, we introduce a class of discrete kernels $p_{n-k}^{(n)}$ as follows
\begin{align}\label{def:DCC kernel}
p_{n-k}^{(n)}:=\sum_{j=k}^n {\theta_{j-k}^{(j)}}
\quad \text{for $1\le k\le n$}.
\end{align}
It follows from \cite[Subsection 2.2]{Liao2020Positive}
that the discrete convolution kernels $p_{n-k}^{(n)}$
are complementary to the original kernels $a_{n-k}^{(n)}$ in the following sense,
\[\sum_{j=k}^np_{n-j}^{(n)}a_{j-k}^{(j)}\equiv 1\quad \text{for $1\le k\le n$}.\]
That is to say, the new kernels $p_{n-k}^{(n)}$ is just
the discrete complementary convolution (DCC) kernels
called by \cite{Liao2018Sharp,Liao2018discrete,Liao2018Unconditional}.

%%%%%%%%%%%%%%%%%%%%%%%%%%%%%%%%%%%%%%%%%%%%%%%%%%%%%%%%%%%%%%%%%%%%%%%%%
\begin{figure}[htb!]
\centering
\includegraphics[width=3.8in]{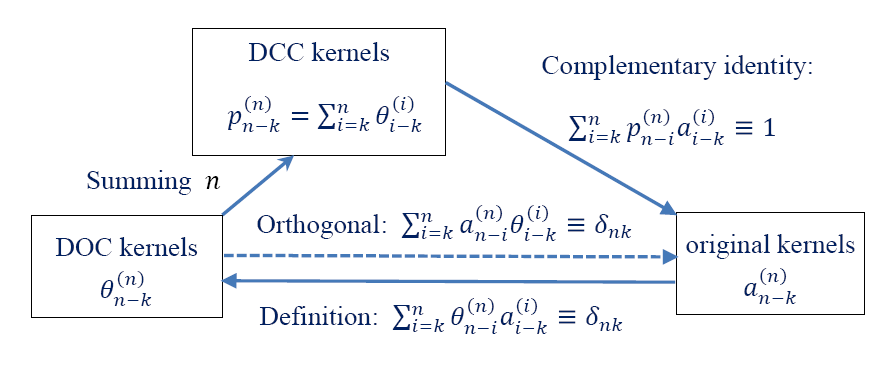}
\caption{The relationship diagram for DOC, DCC and original L1 kernels.}
\label{fig:DOC and DCC relation}
\end{figure}
%%%%%%%%%%%%%%%%%%%%%%%%%%%%%%%%%%%%%%%%%%%%%%%%%%%%%%%%%%%%%%%%%%%%%%%%%

The DOC kernels were originally constructed
in \cite{Liao2019Adaptive} for the numerical analysis of variable-step BDF2
approximation of the first time derivative.
This type of kernels were also applied recently
to handle a general class of discrete kernels in \cite{Liao2020Positive}
and to deal with the discrete L1$_R$ kernels in \cite{Liao2020An}.
The DCC kernels were originally
introduced in \cite{Liao2018Sharp}
to develop the discrete fractional Gr\"onwall inequality.
Figure \ref{fig:DOC and DCC relation} describes
some relationship links between the L1 kernels, DOC and DCC kernels.
In addition, the definition \eqref{def:DCC kernel} gives
the following relationship between the DOC kernels $ \theta_{n-k}^{(n)}$
and the DCC kernels $ p_{n-k}^{(n)}$,
\begin{align}\label{def:equiva DCC kernel}
\theta_{0}^{(n)}={p_{0}^{(n)}}\quad \text{and}\quad
\theta_{n-k}^{(n)}={p_{n-k}^{(n)}}-{p_{n-k-1}^{(n-1)}}
\quad \text{for $1\le k\le n-1$}.
\end{align}
Some properties of the DCC kernels are collected,
see \cite[Lemma 2.2]{Liao2020Positive} and \cite[Lemma 2.5]{Liao2018discrete}.

\begin{lemma}\label{lem:DCC kernel property}
For any $n\ge 2$, the DCC kernels $p_{n-k}^{(n)}$ defined in
\eqref{def:DCC kernel} satisfy
\begin{align*}
 p_{n-k}^{(n)}\ge 0
\quad \text{for $1\le k\le n$;}\quad\text{but}\quad \sum_{j=1}^n {p_{n-j}^{(n)}} \le \omega_{1+\alpha}(t_n).
\end{align*}
%\begin{itemize}
%  \item [(a)] $\displaystyle p_{n-k}^{(n)}=\sum_{j=k}^n {\theta_{j-k}^{(j)}}\ge 0$ for $1\le k\le n$;
%  \item [(b)] $\displaystyle \sum_{j=1}^n {p_{n-j}^{(n)}} \le \omega_{1+\alpha}(t_n).$
%\end{itemize}
\end{lemma}

\begin{lemma}\label{lem:DOC to DCC kernels}
For any real sequence $\{w_k\}_{k=1}^n$,
it holds that
\begin{align*}
2w_n\sum_{k=1}^n{\theta_{n-k}^{(n)}}w_k
\ge \sum_{k=1}^n{p_{n-k}^{(n)}}w_k^2-\sum_{k=1}^{n-1}{p_{n-k-1}^{(n-1)}}w_k^2
+\frac{1}{{\theta_0^{(n)}}}\braB{\sum_{k=1}^n{\theta_{n-k}^{(n)}}w_k}^2
\quad \text{for $ n\ge 1$.}
\end{align*}
\end{lemma}

\begin{proof}
By means of the DOC kernels $\theta_{n-k}^{(n)}$,
we define the following auxiliary kernels
\begin{align}\label{lemp:defini zeta}
\zeta_{n-k}^{(n)}:=\sum_{j=k}^{n}{\theta_{n-j}^{(n)}}\quad\text{for $n\ge1$.}
\end{align}
The sign and sum properties of DOC kernels $\theta_{n-k}^{(n)}$ in Lemma \ref{lem:DOC kernel property} indicate that
\[
\zeta_{n-k}^{(n)}\ge0\quad \text{for $1\le k\le n$}
\quad\text{and}\quad
\zeta_{k-1}^{(n)}\ge\zeta_{k}^{(n)} \quad \text{for $1\le k\le n-1$.}
\]
Consequently, the auxiliary kernels $ \zeta_{n-k}^{(n)} $ satisfy
the assumption of \cite[Lemma A.1]{Liao2018discrete}, which
yields the following inequality
\begin{align}\label{lemp:zeta inequa}
2w_n\sum_{k=1}^{n}\zeta_{n-k}^{(n)}\diff w_k
\ge\sum_{k=1}^{n}\zeta_{n-k}^{(n)}(w_k^2-w_{k-1}^2)
+\frac{1}{\zeta_{0}^{(n)}}\braB{\sum_{k=1}^{n}\zeta_{n-k}^{(n)}\diff w_k}^2.
\end{align}
By taking $w_0=0$,
it is easy to check the following identity
\begin{align*}
\sum_{k=1}^{n}\zeta_{n-k}^{(n)}\diff w_k
=\zeta_{0}^{(n)}w_n
+\sum_{k=1}^{n-1}(\zeta_{n-k}^{(n)}-\zeta_{n-k-1}^{(n)})w_k=\sum_{k=1}^{n}{\theta_{n-k}^{(n)}}w_k,
\end{align*}
where we used the definition \eqref{lemp:defini zeta} in the last step.
Inserting the above identity into \eqref{lemp:zeta inequa} and
applying the relationship \eqref{def:equiva DCC kernel}
to the resulting inequality, we get
the claimed result.
\end{proof}

We are in position to establish a discrete energy law for
the  L1 scheme \eqref{sche:adap AC L1}.
To this end, define a discrete counterpart of
the variational energy \eqref{cont:free energy} as follows
\begin{align*}
\mathcal{E}_\alpha[\phi^0]:=E[\phi^0]\quad\text{and}\quad\mathcal{E}_\alpha[\phi^n]:=E[\phi^n]
+\frac{\kappa}{2}\sum_{j=1}^np_{n-j}^{(n)}\mynormb{\mu^j}^2
\quad\text{for $n\ge 1$,}
\end{align*}
where $E[\phi^n]$ is the discrete version of
the free energy functional \eqref{cont:free energy},
\begin{align*}
E\kbra{\phi^n}:=\frac{\epsilon^2}{2}\mynormb{\nabla_h\phi^n}^2
+\frac14\mynormb{\bra{\phi^n}^2-1}^2\quad\text{for $n\ge 0$.}
\end{align*}

As noted in \cite{Liao2018Sharp,Liao2018discrete}, the DCC kernels $p_{n-j}^{(n)}$
were designed to simulate the continuous kernel
of Riemann-Liouville fractional integral $\mathcal{I}_t^\alpha$.
Comparing the variational energy \eqref{cont:TFAC variational energy} with
the above discrete version $\mathcal{E}_\alpha[\phi^n]$, we see again that
the DCC kernels $p_{n-j}^{(n)}$ ``define"
a certain formula for the Riemann-Liouville fractional integral, that is, $\bra{\mathcal{I}_t^\alpha v}(t_n)\approx\sum_{j=1}^np_{n-j}^{(n)}v^j$.

%\[
%\bra{\mathcal{I}_t^\alpha v}(t_n)\approx\sum_{j=1}^np_{n-j}^{(n)}v^j
%\quad\text{for $n\ge1$.}
%\]
%%%%%%%%%%%%%%%%%%%%%%%%%%%%%%%%%%%%%%%%%%%%%%%%%%%%%%%%%%%%%%%%%%%%%%%
\begin{theorem}\label{thm:energy decay law}
Under the time-step restriction \eqref{restr:solvab time step},
the L1 scheme \eqref{sche:adap AC L1}
preserves the following  discrete energy law at each time level,
\begin{align*}
\partial_{\tau}\mathcal{E}_\alpha[\phi^n]\le0\quad \text{for $1\le n\le N$.}
\end{align*}
\end{theorem}
%%%%%%%%%%%%%%%%%%%%%%%%%%%%%%%%%%%%%%%%%%%%%%%%%%%%%%%%%%%%%%%%%%%%%%%

\begin{proof}
Making the inner product of the first and  second equations of
\eqref{sche:adap AC L1 equi form} by $\mu^n$ and $\diff \phi^n$, respectively,
and adding up the resulting two equalities, one has
\begin{align}\label{thmp:energy law inner}
\kappa\sum_{j=1}^n\theta_{n-j}^{(n)}\myinnerb{\mu^j,\mu^n}
+\myinnerb{\bra{\phi^n}^3-\phi^n,\diff \phi^n}
-\myinnerb{\epsilon^2\Delta_h\phi^n,\diff \phi^n}=0.
\end{align}
An application of the following identity
\begin{align*}
4\brab{a^3-a}\bra{a-b}
=\bra{a^2-1}^2-\bra{b^2-1}^2-2\bra{1-a^2}\bra{a-b}^2+\bra{a^2-b^2}^2
\end{align*}
to the second term of equation \eqref{thmp:energy law inner} yields
\[
\myinnerb{\bra{\phi^n}^3-\phi^n,\diff \phi^n}
\ge \frac14\mynormb{\bra{\phi^n}^2-1}^2-\frac14\mynormb{\bra{\phi^{n-1}}^2-1}^2
-\frac12\mynormb{\diff \phi^n}^2.
\]
Moreover, the discrete Green's formula together with $2a(a-b)=a^2-b^2+(a-b)^2$ gives
\begin{align*}
-\myinnerb{\Delta_h\phi^n,\diff \phi^n}
=\frac12\mynormb{\nabla_h\phi^n}^2
-\frac12\mynormb{\nabla_h\phi^{n-1}}^2
+\frac12\mynormb{\nabla_h\diff \phi^n}^2.
\end{align*}
Inserting the above results into equation \eqref{thmp:energy law inner},
one obtains
\begin{align}\label{thmp:energy law norm}
\kappa\sum_{j=1}^n\theta_{n-j}^{(n)}\myinnerb{\mu^j,\mu^n}
+\frac{\epsilon^2}{2}\mynormb{\nabla_h\diff \phi^n}^2
-\frac12\mynormb{\diff \phi^n}^2
+E\kbra{\phi^n}
\le E\kbra{\phi^{n-1}}.
\end{align}
By means of  Lemma \ref{lem:DOC to DCC kernels},
the first term at the left hand side of  \eqref{thmp:energy law norm}
can be bounded by
\begin{align*}
\kappa\sum_{j=1}^n\theta_{n-j}^{(n)}\myinnerb{\mu^j,\mu^n}
\ge \frac{\kappa}{2}\sum_{j=1}^np_{n-j}^{(n)}\mynormb{\mu^j}^2
-\frac{\kappa}{2}\sum_{j=1}^{n-1}p_{n-1-j}^{(n-1)}\mynormb{\mu^j}^2
+\frac{a_0^{(n)}}{2\kappa}\mynormb{\diff \phi^n}^2,
\end{align*}
where the first equation of \eqref{sche:adap AC L1 equi form} and $\theta_0^{(n)}=1/a_0^{(n)}$
have been used in last term.
%Moreover, the first equation of
%\eqref{sche:adap AC L1 equi form} gives the following identity
%\begin{align*}
%\mynormb{\kappa\sum_{j=1}^n\theta_{n-j}^{(n)}\mu^j}^2
%=\mynormb{\diff \phi^n}^2.
%\end{align*}
Substituting the above results into inequality \eqref{thmp:energy law norm}, one has
\begin{align*}
\braB{\frac{a_0^{(n)}}{2\kappa}-\frac12}\mynormb{\diff \phi^n}^2
+\mathcal{E}_\alpha\kbra{\phi^n}
\le \mathcal{E}_\alpha\kbra{\phi^{n-1}}.
\end{align*}
As a consequence, the time-step condition
\eqref{restr:solvab time step} yields the claimed result immediately.
%This completes the proof.
\end{proof}

%%%%%%%%%%%%%%%%%%%%%%%%%%%%%%%%%%%%%%%%%%%%%%%%%%%%%%%%%%%%%%%%%%%%%%%%

As the fractional order $\alpha\rightarrow 1$,
the definition \eqref{def:DOC kernel}
gives $\theta_0^{(n)}\rightarrow \tau_n$
and $\theta_{n-k}^{(n)}\rightarrow 0$ for $1\le k\le n-1$,
and the definition \eqref{def:DCC kernel} yields
$p_{n-k}^{(n)}\rightarrow \tau_k$ for $1\le k\le n$.
So the discrete energy dissipation law in Theorem \ref{thm:energy decay law} becomes
\[
\partial_{\tau}\mathcal{E}_\alpha[\phi^n]\le0
\longrightarrow
\partial_\tau \bra{E\kbra{\phi}}^n
+\frac{\kappa}{2}\mynormb{\mu^n}^2\le 0
\quad\text{as $\alpha\rightarrow 1$,}
\]
which is just the energy dissipation law \eqref{ieq: energy dissipation law}
 of the backward Euler scheme \eqref{sche: backward Euler AC}
for the classical AC model. In this sense, we say that the energy dissipation law
in Theorem \ref{thm:energy decay law} is asymptotically compatible
 in the fractional order $\alpha\rightarrow1$ limit.
Also, the time-step condition \eqref{restr:solvab time step}
in Theorem \ref{thm:energy decay law} is sharp
since it is asymptotically compatible with the step-size condition
in \eqref{cond: time-step size} as the fractional order $\alpha\rightarrow 1$.

%As a further comment, it would be very interesting to investigate the discrete
%variational energy dissipation law for some other high-order numerical
%formulas of Caputo derivative, such as the Alikhanov formula \cite{Liao2018A},
%the Caputo's BDF2 formula \cite{Liao2016Stability}
%and the recent L1$^+$ formula \cite{Ji2019Adaptive}.

%%%%%%%%%%%%%%%%%%%%%%%%%%%%%%%%%%%%%%%%%%%%%%%%%%%%%%%%%%%%%%%%%%%%%%%%

%%%%%%%%%%%%%%%%%%%%%%%%%%%%%%%%%%%%%%%%%%%%%%%%%%%%%%%%%%%%%%%%%%%%%%%%
\begin{lemma}\label{lem:bound solution}
The solution of the L1 scheme \eqref{sche:adap AC L1} satisfies $\mynormb{\phi^n}_{\ell^6}\le c_0$ for $n\ge1$,
%\[
%\mynormb{\phi^n}_{H_h^1}\le c_0,\quad
%\mynormb{\phi^n}_{\ell^6}\le c_0\quad\text{for $n\ge 1$,}
%\]
where the constant $c_0$ is dependent on the domain $\Omega$,
the parameter $\epsilon$ and the initial value $\phi^0$,
but independent of the time $t_n$, step sizes $\tau_n$
and time-step ratios $r_n$.
\end{lemma}
%%%%%%%%%%%%%%%%%%%%%%%%%%%%%%%%%%%%%%%%%%%%%%%%%%%%%%%%%%%%%%%%%%%%%%%%
\begin{proof}
%The detailed proof is omitted for brevity.
%We only outline the proof of the claimed priori estimate for brevity.
It follows from Theorem \ref{thm:energy decay law} that
$\mynormb{\phi^n}$ and $\mynormb{\nabla_h\phi^n}$
are bounded.
By the discrete Sobolev embedding inequality \cite[Lemma 3.6]{Wang2013Fourth},
\begin{align}\label{inequ:H1 to L6}
\mynormb{v}_{\ell^6}\le c_{\Omega}\mynormb{v}^{\frac13}
\brab{\mynormb{\nabla_hv}+\mynormb{v}}^{\frac23}
\quad\text{for $v\in\mathbb{V}_{h}$,}
\end{align}
one gets the desired estimate immediately
and completes the proof.
\end{proof}

%%%%%%%%%%%%%%%%%%%%%%%%%%%%%%%%%%%%%%%%%%%%%%%%%%%%%%%%%%%%%%%%%%%%%%%%

Under the time step restriction \eqref{restr:solvab time step}, \cite[Theorem 2.1]{Ji2020Simple}
showed that the solution of L1 scheme \eqref{sche:adap AC L1}
preserves the maximum bound principle numerically,
i.e., $\mynormb{\phi^k}_{\infty}\le 1$ for $1\le k\le N$ if
$\mynormb{\phi^0}_{\infty}\le 1$.
The maximum norm error estimate were obtained in \cite[Theorem 3.1]{Ji2020Simple}.
On the other hand, the discrete $L^6$ norm bound in Lemma \ref{lem:bound solution} would be useful to establish
the $L^2$ norm error estimate for other spatial approximations,
such as finite element and spectral methods.

%%%%%%%%%%%%%%%%%%%%%%%%%%%%%%%%%%%%%%%%%%%%%%%%%%%%%%%%%%%%%%%%%%%%%%%%

\section{$L^2$ norm error estimate}	
Under proper assumptions on initial condition such as
$\Phi_0\in H^2(\Omega)\cap H_0^1(\Omega)$,
the TFAC equation \eqref{cont:FAC equ} was proved
to admit a unique solution that fulfills \cite[Theorem 2.2]{Du2019Time},
\[
\mynormb{\Delta^s\partial_t \Phi}_{L^2(\Omega)}\le C_\phi t^{\alpha(1-s)-1}
\quad\text{for $s\in [0,1)$ and $0<t\le T$},
\]
which typically exhibits a singular behavior at an initial time.
Here and hereafter, to facilitate the numerical analysis of
finite difference methods,
we assume that the initial data $\Phi_0$ has the required regularity
and the solution of the TFAC equation \eqref{cont:FAC equ} satisfies
\begin{align}\label{cont:regular sigma AC}
\mynormb{\Phi(t)}_{W^{4,\infty}(\Omega)}\le C_{\phi},\;\; \mynormb{\Phi^{(\ell)}(t)}_{W^{0,\infty}(\Omega)}\le C_{\phi}\brab{1+t^{\sigma-\ell}}
\quad \text{for $0<t\le T$ and $\ell=1,2$,}
\end{align}
where the regularity parameter $\sigma\in(0,1)$ makes our analysis extendable.
It is to mention that such a realistic regularity assumption
on the exact solution of initial-boundary problem with the Caputo time derivative
is standard in numerical analysis \cite{Liao2018Sharp,Liao2018Unconditional,Ji2020Simple,Stynes2017Error,Kopteva2019Error}.
According to the continuous energy law \eqref{cont:TFAC energy law}
and the Sobolev embedding inequality, there exists a positive constant $c_1$
so that $\mynorm{\Phi^n}_{\ell^6}\le c_1$,
which will be used in the convergence analysis below.

\subsection{Global consistency analysis in time}

Denote the local consistency error of the L1 formula \eqref{sche:L1 formula} by
$\Upsilon^{j}:=(\partial_{t}^{\alpha}v)(t_j)-(\partial_\tau^{\alpha}v)^{j}$  for $j\ge 1.$
%\begin{align}\label{schm:local error}
%\Upsilon^{j}:=(\partial_{t}^{\alpha}v)(t_j)-(\partial_\tau^{\alpha}v)^{j}
%\quad\text{for $j\ge 1.$}
%\end{align}
The following result presents a global consistency  error, see \cite[Lemma 3.3, Theorem 3.1]{Liao2018Sharp}.
%%%%%%%%%%%%%%%%%%%%%%%%%%%%%%%%%%%%%%%%%%%%%%%%%%%%%%%%%%%%%%%
\begin{lemma}\label{lem:glabal error estimate}
The global consistency error of the L1 formula \eqref{sche:L1 formula} is bounded by
\begin{align*}
\sum_{j=1}^{n}p_{n-j}^{(n)}\abs{\Upsilon^{j}}
&\le \sum_{k=1}^{n}p_{n-k}^{(n)}a_0^{(k)}G^k
+\sum_{k=1}^{n-1}p_{n-k}^{(n)}a_0^{(k)}G^k\\
&\le C_v\braB{\,\frac{\tau_1^{\sigma}}{\sigma}
+\frac{1}{1-\alpha}\max_{2\le{k}\le{n}}
t_{k}^{\alpha}t_{k-1}^{\sigma-2}\tau_{k}^{2-\alpha}}\quad\text{for $1\le n\le N$,}
\end{align*}
where the local quantity
$G^k :=2\int_{t_{k-1}}^{t_k}\bra{t-t_{k-1}}\abs{v_{tt}}\zd t$ for $1\le k\le n$.
\end{lemma}
%%%%%%%%%%%%%%%%%%%%%%%%%%%%%%%%%%%%%%%%%%%%%%%%%%%%%%%%%%%%%%%

The initial singularity can
be resolved by enforcing the time meshes to
satisfy \cite{Liao2018Sharp,Liao2018Unconditional}:
\begin{enumerate}[itemindent=1em]
\item[\textbf{AG}.] For a mesh parameter $\gamma\geq{1}$,
  there exists mesh-independent constant $C_{\gamma}>0$ such that
  $\tau_k\le \tau \min\{1,C_{\gamma}t_k^{1-1/\gamma}\}$
  for~$1\le k\le N$~and $t_{k}\le C_{\gamma}t_{k-1}$ for $2\le{k}\le{N}$.
\end{enumerate}
The condition \textbf{AG} implies that
the time steps are graded-like near the initial time $t=0$;
while no special structures are imposed on the meshes
when the time is away from $t=0$, except
the maximum time step restriction $\tau_k\le \tau$.
We here note that a typical example of a family of time meshes
satisfying \textbf{AG} is the smoothly graded time grids
$t_k=T(k/N)^{\gamma}$ for $0\le k\le N$.
If \textbf{AG} holds,
we have the following result, see
\cite[Remark 6]{Liao2018Sharp}
and \cite[Lemma 3.3]{Liao2018Unconditional}.

\begin{lemma}\label{lem:glabal error esti graded}
Under the regularity \eqref{cont:regular sigma AC},
if the time meshes satisfy the condition \emph{\textbf{AG}},
then the global consistency error of the L1 formula \eqref{sche:L1 formula}
can be bounded by
\begin{align*}
\sum_{j=1}^{n}p_{n-j}^{(n)}\abs{\Upsilon^{j}}
\le \frac{C_v}{\sigma(1-\alpha)}\tau^{\min\{2-\alpha,\gamma\sigma\}}
\quad\text{for $1\le{n}\le{N}.$}
\end{align*}
\end{lemma}

\subsection{Convergence analysis}
%%%%%%%%%%%%%%%%%%%%%%%%%%%%%%%%%%%%%%%%%%%%%%%%%%%%%%%%%%%%%%%%%%

Let $c_2:=c_{\Omega}\brat{c_0^2+c_0c_1+c_1^2}$ and
$c_3:=2+c_2^{3/2}$, which may be dependent on the domain $\Omega$,
the parameter $\epsilon$ and the initial value $\phi_0$,
but always independent of the time $t_n$, step sizes $\tau_n$ and step ratios $r_n$.
Also, let $r_{*}:=\min_{1\le k\le N}\{1,r_{k}\}$
be the minimum step-ratio.
Recalling the Mittag--Leffler function $E_\alpha(z):=\sum_{k=0}^\infty\frac{z^k}{\Gamma(1+k\alpha)}$,
we have the following convergence result.

\begin{theorem}\label{thm:L2 convergence}
Assume that the solution of \eqref{cont:FAC equ} satisfies
the regular assumption \eqref{cont:regular sigma AC}.
Suppose further that the time-step size restriction
\eqref{restr:solvab time step}
holds such that the adaptive L1 scheme \eqref{sche:adap AC L1}
is unique solvable and energy stable.
If the maximum step size
$\tau\le1/\sqrt[\alpha]{{2\kappa\epsilon^{-1} c_3\Gamma(2-\alpha)}}$,
then the numerical solution $\phi^n$ of the L1 scheme \eqref{sche:adap AC L1}
is convergent in the discrete $L^2$ norm,
\begin{align*}
\mynormb{\Phi^n-\phi^n}\le C_{\phi} E_{\alpha}\brab{2\kappa\epsilon^{-1}c_3t_n^\alpha/r_{*}}
\braB{\frac{\tau_1^\sigma}{\sigma}
+\frac{1}{1-\alpha}\max_{2\le k\le n}
t_k^{\alpha}t_{k-1}^{\sigma-2}\tau_k^{2-\alpha}+h^2} \quad\text{for $1\le{n}\le{N}.$}
\end{align*}
Specially, when the time mesh satisfies \emph{\textbf{AG}}, it holds that
\begin{align*}
\mynormb{\Phi^n-\phi^n}
\le \frac{C_{\phi}}{\sigma(1-\alpha)}E_{\alpha}
\brab{2\kappa\epsilon^{-1}c_3t_n^\alpha/r_{*}}
\braB{\tau^{\min\{2-\alpha,\gamma\sigma\}}+h^2}
 \quad\text{for $1\le{n}\le{N},$}
\end{align*}
and the  optimal accuracy $O(\tau^{2-\alpha})$ can be achieved
when $\gamma\geq\max{\{1,\,(2-\alpha)/\sigma\}}$.
\end{theorem}
%%%%%%%%%%%%%%%%%%%%%%%%%%%%%%%%%%%%%%%%%%%%%%%%%%%%%%%%%%%%%%%%%%

\begin{proof}
Let the error function
$e_h^n:=\Phi_h^n-\phi_h^n \in \mathbb{V}_{h} $.
Substituting the exact solution $\Phi_h^n$ into the numerical
scheme \eqref{sche:adap AC L1}, one has
\begin{align}\label{thmp:exact equ}
(\partial^{\alpha}_\tau \Phi_h)^n
=\kappa\kbrab{\epsilon^2\Delta_h \Phi_h^n-f(\Phi_h^n)}
+\Upsilon_h^n+\xi_h^n\quad\text{for $n\ge 1$,}
\end{align}
with the initial data $\phi_h^0=\Phi_0(\mathbf{x}_h)$.
Here, the notations $\Upsilon_h^n$ and $\xi_h^n$ represent
the temporal and spatial truncation errors, respectively.
Subtracting the numerical scheme \eqref{sche:adap AC L1}
from the equation \eqref{thmp:exact equ},
we obtain the following error system
\begin{align}\label{thmp:error equ}
(\partial^{\alpha}_{\tau}e_h)^n
=\kappa\bra{\epsilon^2\Delta_he_h^n-f_\phi^ne_h^n}
+\Upsilon_h^n+\xi_h^n\quad\text{for $n\ge 1$,}
\end{align}
where the nonlinear term $f_\phi^n$ is defined by
\[
f_\phi^n:=(\Phi_h^n)^2+\Phi_h^n\phi_h^n+(\phi_h^n)^2-1.
\]
Taking the inner product of \eqref{thmp:error equ} by $e^n$, one obtains
\begin{align}\label{thmp:error equa inner}
\myinnerb{(\partial^{\alpha}_{\tau}e)^n,e^n}
+\kappa\epsilon^2\mynormb{\nabla_he^n}^2
=-\kappa\myinnerb{f_\phi^ne^n,e^n}
+\myinnerb{\Upsilon^{n}+\xi^{n},e^n}\quad\text{for $n\ge 1$,}
\end{align}
where the discrete Green's formula has been used in the above derivation.
Thanks to Lemma \ref{lem:L1 kernel property},
the first term at the left hand side of \eqref{thmp:error equa inner}
can be bounded by \cite{Liao2018Sharp,Liao2018Unconditional},
\begin{align*}
\myinnerb{\bra{\partial_\tau^\alpha e}^n,e^n}
=\sum_{k=1}^{n}a_{n-k}^{(n)}\myinnerb{\diff e^k,e^n}
\ge\mynormb{e^n}\sum_{k=1}^{n}a_{n-k}^{(n)}\diff \mynormb{e^k}.
\end{align*}
For the nonlinear term of the right hand side of equation \eqref{thmp:error equa inner},
one can apply the H\"older inequality
to derive the following estimation
\begin{align*}
\mynormb{f_{\phi}^ne^n}
&\le\mynormb{e^n}+\bra{\mynormb{\Phi^n}_{\ell^6}^2
+\mynormb{\Phi^n}_{\ell^6}\mynormb{\phi^n}_{\ell^6}
+\mynormb{\phi^n}_{\ell^6}^2}\mynormb{e^n}_{\ell^6}\\
&\le \mynormb{e^n}+c_2\mynormb{e^n}^{\frac13}
\braB{\mynormb{\nabla_he^n}+\mynormb{e^n}}^{\frac23},
\end{align*}
where the discrete Sobolev embedding inequality \eqref{inequ:H1 to L6}
has been used.
Then one can use the above estimate and the Young's inequality
to  derive that
\begin{align*}
\abs{\myinnerb{f_\phi^ne^n,e^n}}
&\le \mynormb{e^n}^2+c_2\mynormb{e^n}^{\frac43}
\braB{\mynormb{\nabla_he^n}+\mynormb{e^n}}^{\frac23}\\
&\le \mynormb{e^n}^2+ c_2^{3/2}\epsilon^{-1}\mynormb{e^n}^2
+\epsilon^2\bra{\mynormb{\nabla_he^n}^2+\mynormb{e^n}^2}\\
&\le c_3\epsilon^{-1}\mynormb{e^n}^2
+\epsilon^2\mynormb{\nabla_he^n}^2,
\end{align*}
where the width of diffusive interface $\epsilon\ll 1$
has been used in the last step.
Collecting the above estimates into \eqref{thmp:error equa inner}, one gets
\begin{align*}%\label{thmp:error norm esti}
\mynormb{e^n}\sum_{k=1}^na_{n-k}^{(n)}\diff\mynormb{e^k}
\le \kappa c_3\epsilon^{-1}\mynormb{e^n}^2
+\bra{\mynormb{\Upsilon^n}+\mynormb{\xi^n}}\mynormb{e^n},
\end{align*}
which implies the following inequality
\begin{align*}
\sum_{k=1}^na_{n-k}^{(n)}\diff\mynormb{e^k}
\le \kappa c_3\epsilon^{-1}\mynormb{e^n}
+\mynormb{\Upsilon^n}+\mynormb{\xi^n}.
\end{align*}
If $\tau\le1/\sqrt[\alpha]{{2\kappa\epsilon^{-1} c_3\Gamma(2-\alpha)}}$,
the discrete fractional Gr\"{o}nwall inequality
\cite[Theorem 3.2]{Liao2018discrete}
with the substitutions $\lambda:=\kappa c_3\epsilon^{-1}$,
$v^k:=\mynormb{e^k}$ and $g^n:=\mynormb{\Upsilon^n}+\mynormb{\xi^n}$ yields
\begin{align*}
\mynormb{e^n}\le 2E_{\alpha}\brab{2\kappa\epsilon^{-1}c_3t_n^\alpha/r_{*}}
\braB{\mynormb{e^0}+\max_{1\le k\le n}\sum_{j=1}^k{p_{k-j}^{(k)}}
\mynormb{\Upsilon^j}
+C_\phi\omega_{1+\alpha}(t_{n})h^2},
\end{align*}
where the regularity condition \eqref{cont:regular sigma AC}
has been used to derive that $\mynormb{\xi^n}\le C_\phi h^2$.
The claimed error estimate follows from Lemmas \ref{lem:glabal error estimate}
and \ref{lem:glabal error esti graded} immediately. %This completes the proof.
\end{proof}

\section{Numerical examples}

In this section, we present several numerical examples
to illustrate the efficiency and accuracy of the adaptive
L1 method \eqref{sche:adap AC L1} for the TFAC equation \eqref{cont:FAC equ}.
At each time level, the nonlinear scheme is solved by employing a simple fixed-point algorithm
with the termination error $10^{-12}$.
In addition,
the sum-of-exponentials technique \cite{Jiang2017Fast}
with the absolute tolerance error $\epsilon=10^{-12}$
and cut-off time $\Delta{t}=\tau_1$
is always used to speed up the evaluation
of the L1 formula \eqref{sche:L1 formula}.

\subsection{Accuracy verification}

\begin{example}\label{examp:accuracy test}
Consider an exact solution
$\Phi(\mathbf{x},t)=\omega_{1+\sigma}(t)\sin({x})\sin({y})$
with $\sigma\in(0,1)$ by adding an exterior force to the TFAC equation \eqref{cont:FAC equ}.
\end{example}

%%%%%%%%%%%%%%%%%%%%%%%%%%%%%%%%%%%%%%%%%%%%%%%%%%%%%%%%%%%%%%%%%%%
\begin{table}[htb!]
\begin{center}
\caption{Numerical accuracy of L1 scheme \eqref{sche:adap AC L1} for $\alpha=0.4,\,\sigma=0.4$ with $\gamma_{\mathrm{opt}}=4$}\label{examp:error alph 04} \vspace*{0.3pt}
\def\temptablewidth{1.0\textwidth}
{\rule{\temptablewidth}{0.5pt}}
\begin{tabular*}{\temptablewidth}{@{\extracolsep{\fill}}cccccccccc}
\multirow{2}{*}{$N$} &\multirow{2}{*}{$\tau$} &\multicolumn{2}{c}{$\gamma=3 $} &\multirow{2}{*}{$\tau$} &\multicolumn{2}{c}{$\gamma=4 $} &\multirow{2}{*}{$\tau$}&\multicolumn{2}{c}{$\gamma=5 $} \\
             \cline{3-4}          \cline{6-7}         \cline{9-10}
         &         &$e(N)$   &Order &         &$e(N)$   &Order&         &$e(N)$   &Order\\
\midrule
  40     &6.17e-02 &5.03e-02 &$-$   &6.68e-02 &1.35e-02 &$-$  &6.48e-02 &7.52e-03 &$-$\\
  80     &3.07e-02 &2.19e-02 &1.19  &3.35e-02 &4.44e-03 &1.61 &3.28e-02 &2.87e-03 &1.42\\
 160     &1.48e-02 &9.54e-03 &1.14  &1.59e-02 &1.47e-03 &1.49 &1.69e-02 &9.30e-04 &1.69\\
 320     &7.79e-03 &4.15e-03 &1.30  &7.89e-03 &4.88e-04 &1.56 &8.53e-03 &3.16e-04 &1.58\\
 \midrule
\multicolumn{3}{l}{$\min\{\gamma\sigma,2-\alpha\}$}   &1.20 & & &1.60 & &         &1.60\\
\end{tabular*}
{\rule{\temptablewidth}{0.5pt}}
\end{center}
\end{table}
%%%%%%%%%%%%%%%%%%%%%%%%%%%%%%%%%%%%%%%%%%%%%%%%%%%%%%%%%%%%%%%%%%%

%%%%%%%%%%%%%%%%%%%%%%%%%%%%%%%%%%%%%%%%%%%%%%%%%%%%%%%%%%%%%%%%%%%
\begin{table}[htb!]
\begin{center}
\caption{Numerical accuracy of L1 scheme \eqref{sche:adap AC L1} for $\alpha=0.8,\,\sigma=0.4$ with $\gamma_{\mathrm{opt}}=3$}\label{examp:error alph 08} \vspace*{0.3pt}
\def\temptablewidth{1.0\textwidth}
{\rule{\temptablewidth}{0.5pt}}
\begin{tabular*}{\temptablewidth}{@{\extracolsep{\fill}}cccccccccc}
\multirow{2}{*}{$N$} &\multirow{2}{*}{$\tau$} &\multicolumn{2}{c}{$\gamma=2 $} &\multirow{2}{*}{$\tau$} &\multicolumn{2}{c}{$\gamma=3 $} &\multirow{2}{*}{$\tau$}&\multicolumn{2}{c}{$\gamma=4 $} \\
             \cline{3-4}          \cline{6-7}         \cline{9-10}
     &         &$e(N)$   &Order &         &$e(N)$   &Order &         &$e(N)$   &Order\\
\midrule
  40 &5.12e-02 &1.92e-01 &$-$   &6.37e-02 &7.35e-02 &$-$   &6.67e-02 &6.20e-02 &$-$\\
  80 &2.83e-02 &1.10e-01 &0.94  &3.03e-02 &3.30e-02 &1.08  &3.12e-02 &2.75e-02 &1.07\\
 160 &1.41e-02 &6.39e-02 &0.79  &1.61e-02 &1.46e-02 &1.29  &1.66e-02 &1.20e-02 &1.31\\
 320 &7.03e-03 &3.67e-02 &0.80  &7.86e-03 &6.42e-03 &1.14  &8.06e-03 &5.30e-03 &1.13\\
 \midrule
\multicolumn{3}{l}{$\min\{\gamma\sigma,2-\alpha\}$}   &0.8 & & &1.20 & &       &1.20\\
\end{tabular*}
{\rule{\temptablewidth}{0.5pt}}
\end{center}
\end{table}
%%%%%%%%%%%%%%%%%%%%%%%%%%%%%%%%%%%%%%%%%%%%%%%%%%%%%%%%%%%%%%%%%%%

We use $512^2$ spatial meshes to discretize the domain $(0,2\pi)^2$,
 and take $\kappa=1$, $\epsilon^2=0.5$ and $T=1$.
The time interval $[0,T]$ is always divided into
two parts $[0, T_0]$ and $[T_0, T]$ with total $N$ subintervals,
where $T_0=\min\{1/\gamma,T\}$ and $N_0=\lceil \frac{N}{T+1-\gamma^{-1}}\rceil$.
The graded meshes $t_{k}=T_{0}(k/N_0)^{\gamma}$ are applied to
the initial part $[0,T_{0}]$ for resolving the initial singularity.
Let $N_1:=N-N_0$, $S_1=\sum_{k=1}^{N_1}\epsilon_{k}$
and $\epsilon_{k}\in(0,1)$ be the random numbers.
The random time-steps
$\tau_{N_{0}+k}:=(T-T_{0})\epsilon_{k}/S_1$ for $1\le k\le N_1$
are uniformly distributed in the remainder interval $[T_{0},T]$
to test the mesh-robustness of L1 scheme.

We only test the time accuracy.
The discrete $L^2$ errors
$e(N):=\max_{1\le n\le N}\mynormb{\Phi^n-\phi^n}$
and the experimental order of convergence is estimated by
$$\text{Order}\approx\frac{\log\bra{e(N)/e(2N)}}{\log\bra{\tau(N)/\tau(2N)}},$$
where $\tau(N)$ denotes the maximum time-step size for total $N$ subintervals.
The numerical results in Tables \ref{examp:error alph 04}
and \ref{examp:error alph 08} are computed for $\sigma=\alpha=0.4$
and $\sigma=0.4,\,\alpha=0.8$, respectively, with different grading parameters $\gamma$.
As observed, the L1 scheme \eqref{sche:adap AC L1} is of order $O(\tau^{\gamma \sigma})$
if the graded parameters $\gamma<\gamma_{\text{opt}}$,
while the optimal accuracy $ O(\tau^{2-\alpha})$ can be reached
when $\gamma\ge\gamma_{\text{opt}}$.
Experimentally, they support the sharpness of our theoretical findings.

\subsection{Simulation of coarsening dynamics}
\begin{example}
The coarsening dynamics of
the TFAC model is examined with  $\kappa=1$ and
$\epsilon=0.05$.
The initial condition is generated as
$\Phi_{0}(\mathbf{x})=\mathrm{rand}(\mathbf{x}) $,
where $\mathrm{rand}(\mathbf{x}) $ is uniformly distributed random number
varying from $ -0.001 $ to 0.001 at each grid points.
\end{example}

%%%%%%%%%%%%%%%%%%%%%%%%%%%%%%%%%%%%%%%%%%%%%%%%%%%%%%%%%%%%%%%%%%
\begin{figure}[htb!]
\centering
\includegraphics[width=1.47in]{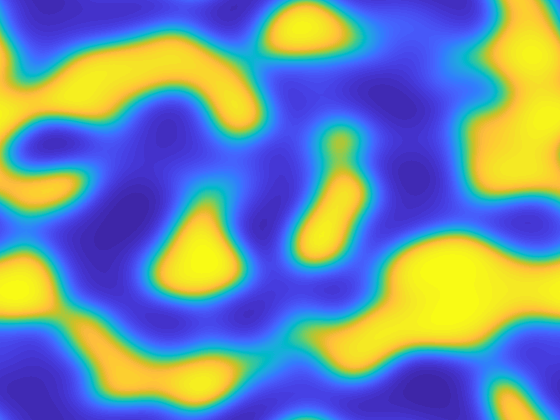}
\includegraphics[width=1.47in]{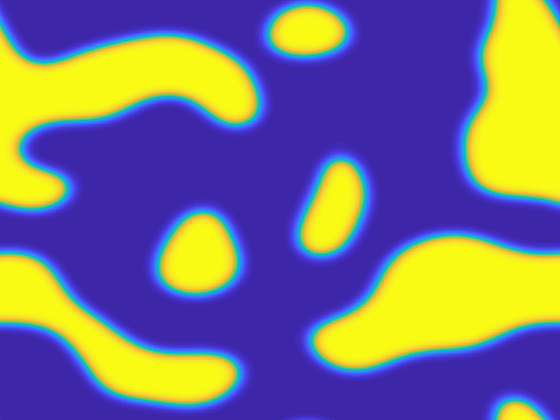}
\includegraphics[width=1.47in]{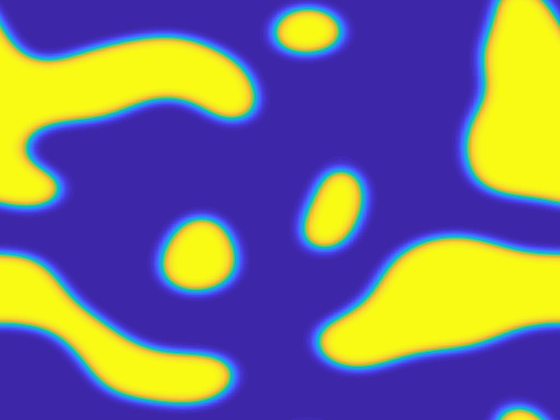}
\includegraphics[width=1.47in]{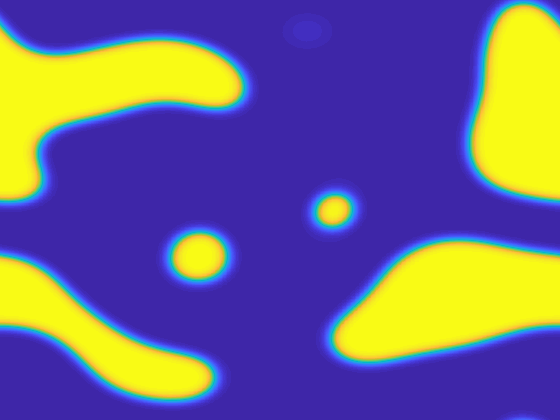}\\
\includegraphics[width=1.47in]{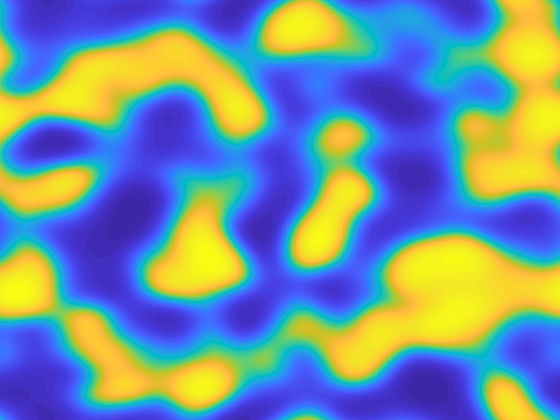}
\includegraphics[width=1.47in]{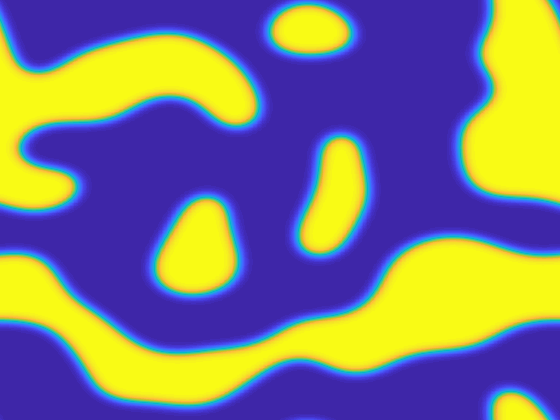}
\includegraphics[width=1.47in]{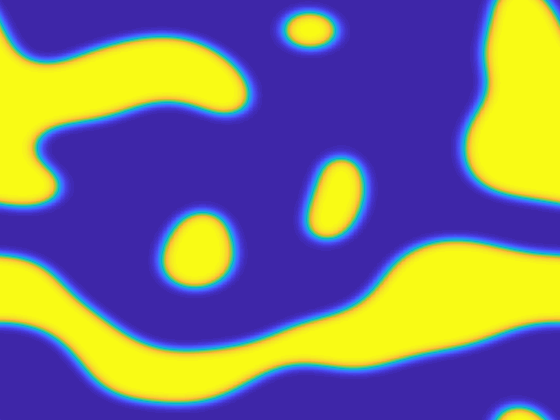}
\includegraphics[width=1.47in]{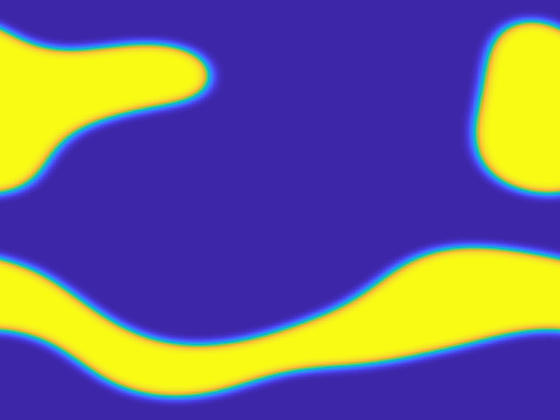}\\
\includegraphics[width=1.47in]{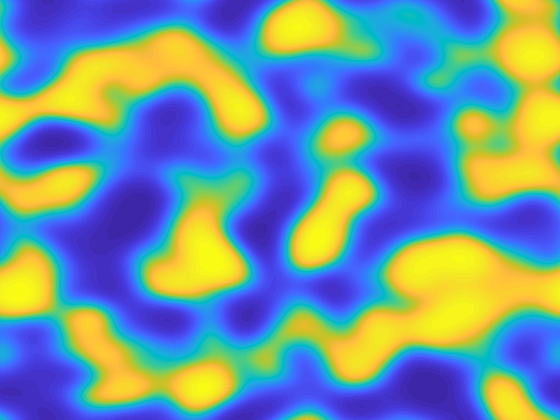}
\includegraphics[width=1.47in]{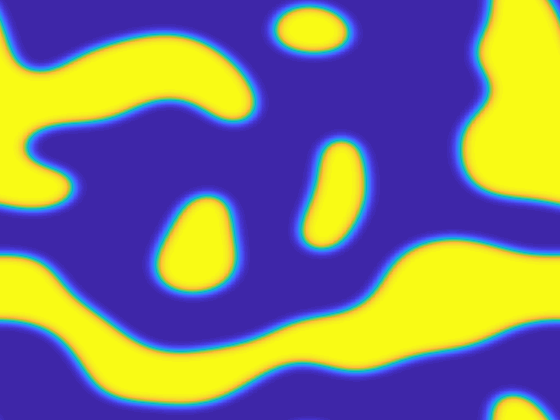}
\includegraphics[width=1.47in]{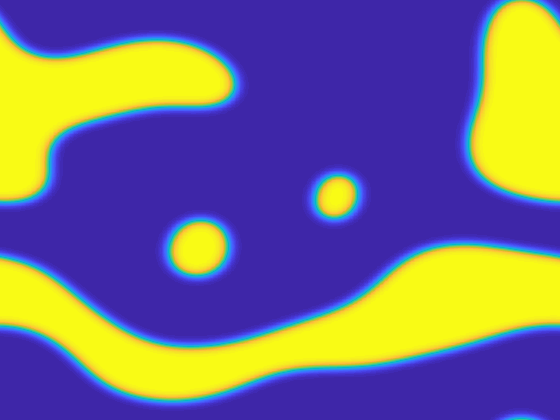}
\includegraphics[width=1.47in]{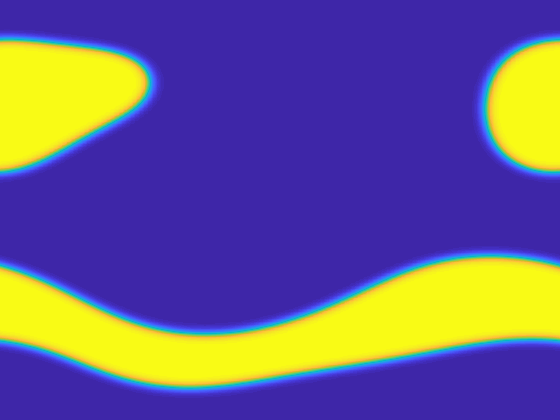}
\caption{Profiles of coarsening dynamics at $t=10, 50, 100, 300$ (from left to right) with different
  fractional orders $\alpha=0.4, 0.7, 0.9$ (from top to bottom), respectively.}
  \label{exam:AC dynamic snap}
\label{figuer}
\end{figure}
%%%%%%%%%%%%%%%%%%%%%%%%%%%%%%%%%%%%%%%%%%%%%%%%%%%%%%%%%%%%%%%%%%%

%%%%%%%%%%%%%%%%%%%%%%%%%%%%%%%%%%%%%%%%%%%%%%%%%%%%%%%%%%%%%%%%%%%
\begin{figure}[htb!]
\centering
\includegraphics[width=2.0in]{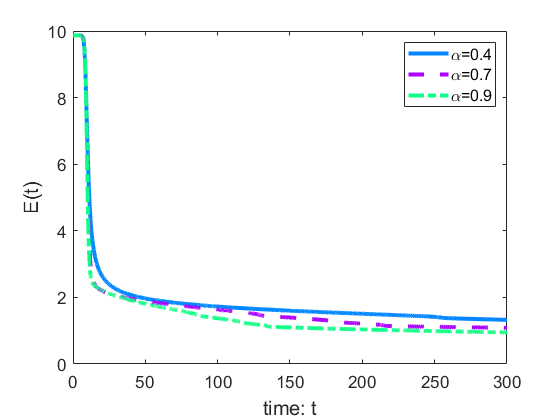}
\includegraphics[width=2.0in]{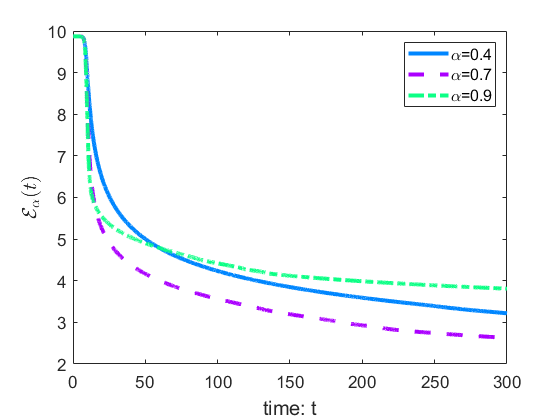}
\includegraphics[width=2.0in]{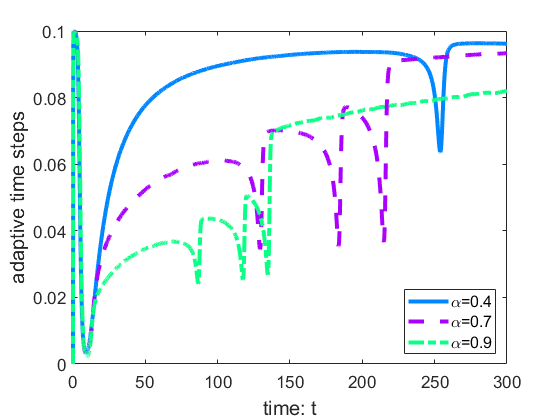}
\caption{$E(t)$, $\mathcal{E}_\alpha(t)$
 and adaptive time steps in the coarsening dynamics
 with the difference fractional orders $\alpha=0.4, 0.7, 0.9$ (from top to bottom), respectively.}
\label{exam:AC dynamic energy}
\end{figure}
%%%%%%%%%%%%%%%%%%%%%%%%%%%%%%%%%%%%%%%%%%%%%%%%%%%%%%%%%%%%%%%%%%%

Several practical implementations related numerical simulations are listed below.
We use $128^2$ grids in space to discretize
the domain $(0,2\pi)^2$. The graded mesh $t_k=T_0(k/N_0)^{\gamma}$ with
the settings $\gamma=3,N_0=30$ and $T_0=0.01$
are always applied in the initial part $[0,T_0]$,
cf. \cite{Ji2020Simple,Ji2019Adaptive}.
For the remainder interval,
the time steps are selected according to
the following adaptive time-stepping strategy
\cite{Zhang2012An,Huang2020Parallel}
\begin{align*}
\tau_{ada}
=\max\Bigg\{\tau_{\min},
\frac{\tau_{\max}}{\sqrt{1+\eta\mynormb{\partial_\tau \phi^n}^2}}\Bigg\}.
\end{align*}
Always, we take
$\tau_{\max}=10^{-1},\,\tau_{\min}=10^{-3}$ and the user parameter $\eta=10^3$.

The  profiles of coarsening dynamics with different
fractional orders $\alpha=0.4,0.7$ and 0.9 for the TFAC model are shown
in Figure \ref{exam:AC dynamic snap}.
Snapshots are taken at time $t=10, 50, 100$ and 300, respectively.
We observe that the coarsening rates are dependent on the fractional order and the time period.
Near the initial time, the small the fractional order $\alpha$, the faster
the coarsening dynamics; while the time goes away from the initial time,
the small the fractional order $\alpha$, the slower
the coarsening dynamics.
The calculated energies using the adaptive time steps for
the coarsening dynamics in Figure \ref{exam:AC dynamic snap} are depicted
in Figure \ref{exam:AC dynamic energy}.
We observe that both the original and variational energies decay with respect to time,
while the former decays faster for larger $\alpha$.

\section*{Acknowledgements}
The authors would like to thank Dr. Bingquan Ji for his help on numerical computations.

%%%%%%%%%%%%%%%%%%%%%%%%%%%%%%%%%%%%%%%%%%%%%%%%%%%%%%%%%%%%%%%%%%%%%%%
%\bibliographystyle{unsrt}%ÉèÖÃ²Î¿¼ÎÄÏ×ÑùÊ½
\bibliographystyle{plain}%ÉèÖÃ²Î¿¼ÎÄÏ×ÑùÊ½
\bibliography{FAC-Energy-StableL1}

\end{document}